\documentclass[12pt]{iopart}

\usepackage{amssymb}

\newtheorem{thm}{Theorem}[section]
\newtheorem{rk}[thm]{Remark} 
\newtheorem{prop}[thm]{Proposition} 
\newtheorem{clly}[thm]{Corollary} 
\newtheorem{lemma}[thm]{Lemma}
 
\newtheorem{defi}[thm]{Definition} 
 
\newtheorem{ex}[thm]{Example}

\begin{document}

\title[Expansivity of ergodic measures with positive entropy]{Expansivity of ergodic measures
with positive entropy}

\author{A. Arbieto, C. A. Morales}

\address{Instituto de Matem\'atica, Universidade Federal do Rio de Janeiro, P. O.
Box 68530, 21945-970 Rio de Janeiro, Brazil}
\ead{arbieto@im.ufrj.br,morales@impa.br}

\begin{abstract}
We prove that for every ergodic invariant measure with positive
entropy of a continuous map on a compact metric space
there is $\delta>0$ such that the dynamical $\delta$-balls have measure zero.
We use this property to prove, for instance, that
the stable classes have measure zero with respect to any ergodic invariant measure with positive entropy.
Moreover, continuous maps which either have countably many stable classes or
are Lyapunov stable on their recurrent sets
have zero topological entropy.
We also apply our results to the Li-Yorke chaos.
\end{abstract}


\vspace{2pc}


\section{Introduction}

\noindent
Ergodic measures with positive entropy for continuous maps on compact metric spaces
have been studied in the recent literature.
For instance,
\cite{bhr} proved that
the set of points belonging to a proper asymptotic pair
(i.e. points whose stable classes are not singleton)
constitute a full measure set.
Moreover, \cite{h} proved
that if $f$ is a homeomorphism with positive entropy $h_\mu(f)$
with respect to one of such measures $\mu$, then
there is a full measure set $A$ such that for all $x\in A$
there is a closed subset $A(x)$ in the stable class of $x$
satisfying $h(f^{-1},A(x))\geq h_\mu(f)$,
where $h(\cdot,\cdot)$ is the Bowen's entropy operation \cite{b}.
We can also mention
\cite{cj} which proved that every
ergodic endomorphism on a Lebesgue probability space having positive entropy on
finite measurable partitions formed by continuity sets
is pairwise sensitive.

In this paper we prove that
these measures have an additional property closely related to \cite{cj}.
More precisely, we prove that
for every ergodic invariant measure with positive
entropy of a continuous map on a compact metric space
there is $\delta>0$ such that the dynamical $\delta$-balls have measure zero.
Measures with this property will be called {\em expansive measures}
in part motivated by the classical definition of {\em expansive map} which requires $\delta>0$ such that
{\em the dynamical $\delta$-balls reduce to singleton} \cite{u}.

We shall prove that expansive measures
satisfy certain properties which are interesting by themselves.
With the aid of these properties we prove that, on compact metric spaces,
every stable class has measure zero with respect to any ergodic measure with positive entropy
(this seems to be new as far as we know).
We also prove through the use of expansive measures that
every continuous map on a compact metric space
exhibiting countably many stable classes
has zero topological entropy
(a similar result with different techniques has been obtained in \cite{hy} but in the {\em transitive} case).
Still in the compact case
we prove that
every continuous map which is Lyapunov stable on its recurrent set
has zero topological entropy too
(this result is well-known but for {\em one-dimensional maps} \cite{fs}, \cite{si}, \cite{z}).
Finally we use expansive measures to give necessary conditions for a
continuous map on a Polish space to be chaotic in the sense of Li and Yorke \cite{liy}.

\section{Definition of expansive measure}

\noindent
In this section we define expansive measures,
present some examples and
give necessary and sufficient
conditions for a Borel probability to be expansive.

To motivate the definition we recall the notion of expansive map first.
Let $(X,d)$ a metric space and $f:X \to X$ be a map of $X$.
Given $x\in X$ and $\delta>0$ we define the dynamical ball of radio $\delta$,
$$
\Phi_\delta(x)=\{y\in X:d(f^i(y),f^i(x))\leq \delta, \quad\forall i\in I\!\! N\}.
$$
We say that $f$ is {\em expansive} (or {\em positively expansive})
if there is $\delta>0$ (called {\em expansivity constant}) such that
for all $x\neq y$ in $X$ there is $n\in\mathbb{N}$ such that $d(f^n(x),f^n(y))> \delta$
(c.f. \cite{hk}). Equivalently, $f$ is expansive if there is $\delta>0$ such that
$$
\Phi_\delta(x)=\{x\},
\quad \quad \forall x\in X.
$$
In such a case every nonatomic Borel measure $\mu$ of $X$ satisfies
\begin{equation}
\label{paratodo}
\mu(\Phi_\delta(x))=0,\quad\quad
\forall x\in X.
\end{equation}
This suggests the following definition in which the term
{\em measurable} means measurable with respect to the Borel $\sigma$-algebra.

\begin{defi}
\label{expansive-measure}
An {\em expansive measure}(\footnote{the name {\em positively expansive measure} seems to be the correct one.
Another possible names are
{\em pairwise sensitive measure} or
{\em symmetrically sensitive measure}.}) of a measurable map $f:X\to X$ is a
Borel measure $\mu$ for which there is $\delta>0$ satisfying
(\ref{paratodo}).
\end{defi}

Notice that this definition does not assume that
the map $f$ (resp. the measure $\mu$) is {\em measure-preserving}
(resp. {\em invariant}), i.e., $\mu=\mu\circ f^{-1}$.
In fact, this hypothesis will not be assumed unless otherwise stated.
Here are some examples.

\begin{ex}
\label{exa1}
Neither the identity nor the constant maps on separable metric spaces
have expansive measures.
On the other hand, every expansive measure is nonatomic.
The converse is true for expansive maps, i.e.,
every nonatomic Borel probability measure is expansive with respect to
any expansive map.
On the other hand, there are nonexpansive continuous maps on certain compact metric spaces for
which every nonatomic measure is expansive \cite{m1}.
The homeomorphism $f(x)=2x$ in $\mathbb{R}$ exhibits expansive measures
(e.g. the Lebesgue measure) but
not expansive {\em invariant} ones.
\end{ex}

Hereafter all expansive measures will be probability measures.
Now we present an useful characterization of expansive measures.

\begin{lemma}
 \label{forall}
A Borel probability measure $\mu$ is expansive for a measurable map $f$
if and only if there is $\delta>0$ such that
\begin{equation}
\label{quasitodo}
\mu(\Phi_\delta(x))=0,
\quad\quad
\forall \mu\mbox{-a.e. }
x\in X.
\end{equation}
\end{lemma}

\noindent{{\bf Proof.}}
We only have to prove that (\ref{quasitodo}) implies that $\mu$ is expansive.
Fix $\delta>0$ satisfying (\ref{quasitodo})
and suppose by contradiction $\mu$ is not expansive.
Then,
there is $x_0\in X$ such that $\mu(\Phi_{\delta/2}(x_0))>0$.
Denote $X_\delta=\{x\in X:\mu(\Phi_\delta(x))=0\}$
so $\mu(X_\delta)=1$.
Since $\mu$ is a probability we obtain
$X_\delta\cap \Phi_{\frac{\delta}{2}}(x_0)\neq\emptyset$ so
there is $y_0\in \Phi_{\frac{\delta}{2}}(x_0)$ such that $\mu(\Phi_\delta(y_0))=0$.
Now if $x\in \Phi_{\frac{\delta}{2}}(x_0)$ we have $d(f^i(x),f^i(x_0))\leq \frac{\delta}{2}$ (for all $i\in \mathbb{N}$) and, since
$y_0\in \Phi_{\frac{\delta}{2}}(x_0)$, we obtain
$d(f^i(y_0),f^i(x_0))\leq \frac{\delta}{2}$ (for all $i\in \mathbb{N}$)
so
$d(f^i(x),f^i(y_0))\leq d(f^i(x),f^i(x_0))+d(f^i(x_0),f^i(y_0))\leq \frac{\delta}{2}+\frac{\delta}{2}=\delta$
(for all $i\in \mathbb{N}$)
proving $x\in \Phi_{\frac{\delta}{2}}(y_0)$.
Therefore $\Phi_{\frac{\delta}{2}}(x_0)\subset \Phi_\delta(y_0)$
so $\mu(\Phi_{\frac{\delta}{2}}(x_0))\leq \mu(\Phi_\delta(y_0))=0$ which is absurd.
This proves the result.
\opensquare

\vspace{5pt}

This lemma together with the corresponding definition for expansive maps suggests the following.

\begin{defi}
\label{sensitivity-constant}
An {\em expansivity constant} of an expansive measure $\mu$ is a constant $\delta>0$
satisfying either (\ref{paratodo}) or (\ref{quasitodo}).
\end{defi}

\section{Some properties of expansive measures}

\noindent
In this section we select the properties of expansive measures
we shall use in the next section.
For the first one we need the following definition.

\begin{defi}
 \label{def-stable-set}
Given a map $f: X\to X$ and
$p\in X$ we define $W^s(p)$, the {\em stable set} of $p$, as the set of points $x$ for which the
pair $(p,x)$ is asymptotic, i.e.,
$$
W^s(p)=\left\{x\in X:\lim_{n\to\infty}d(f^n(x),f^n(p))=0\right\}.
$$
By a {\em stable class} we mean a subset equals to $W^s(p)$ for some $p\in X$.
\end{defi}

The following shows that every stable class is negligible
with respect to any expansive {\em invariant} measure.

\begin{prop}
\label{asymptotic-pair}
The stable classes of a measurable map
have measure zero with respect to any expansive {\em invariant} measure.
\end{prop}

\noindent{{\bf Proof.}}
Let $f:X\to X$ a measurable map and $\mu$ be an expansive invariant measure.
Denoting by $B[\cdot,\cdot]$ the closed ball operation one gets
$$
W^s(p)=\bigcap_{i\in \mathbb{N}^+}\bigcup_{j\in \mathbb{N}}\bigcap_{k\geq j}f^{-k}\left(
B\left[f^k(p),\frac{1}{i}\right]\right).
$$
As clearly
$$
\bigcup_{j\in \mathbb{N}}\bigcap_{k\geq j}f^{-k}\left(B\left[f^k(p),\frac{1}{i+1}\right]\right)
\subseteq \bigcup_{j\in \mathbb{N}}\bigcap_{k\geq j}f^{-k}\left(B\left[f^k(p),\frac{1}{i}\right]\right),
\quad\forall i\in\mathbb{N}^+,
$$
we obtain
\begin{equation}
\label{zero}
\mu(W^s(p))\leq
\lim_{i\to\infty}
\sum_{j\in \mathbb{N}}\mu\left(\bigcap_{k\geq j}f^{-k}\left(B\left[f^k(p),\frac{1}{i}\right]\right)
\right).
\end{equation}
On the other hand,
$$
\bigcap_{k\geq j}f^{-k}\left(B\left[f^k(p),\frac{1}{i}\right]\right)=f^{-j}\left(\Phi_{\frac{1}{i}}(f^j(p))\right)
$$
so
$$
\mu\left(
\bigcap_{k\geq j}f^{-k}\left(B\left[f^k(p),\frac{1}{i}\right]\right)
\right)=\mu\left(
f^{-j}\left(\Phi_{\frac{1}{i}}(f^j(p))\right)
\right)=\mu\left(\Phi_{\frac{1}{i}}(f^j(p))\right)
$$
since $\mu$ is invariant. Then,
taking $i$ large, namely, $i>\frac{1}{\epsilon}$
where $\epsilon$ is a expansivity constant of $\mu$ (c.f. Definition \ref{sensitivity-constant})
we obtain
$\mu\left(\Phi_{\frac{1}{i}}(f^j(p))\right)=0$
so
$$
\mu\left(
\bigcap_{k\geq j}f^{-k}\left(B\left[f^k(p),\frac{1}{i}\right]\right)
\right)=0.
$$
Replacing in (\ref{zero}) we get the result.
\opensquare

\vspace{5pt}

For the second property we will use the following definition \cite{fs}.

\begin{defi}
A map $f: X\to X$ is said to be {\em Lyapunov stable} on $A\subset X$
if for any $x\in A$ and any $\epsilon>0$
there is a neighborhood $U(x)$ of $x$ such that
$d(f^n(x),f^n(y))< \epsilon$ whenever $n\geq 0$ and $y\in U(x)\cap A$.
\end{defi}

(Notice the difference between this definition and the corresponding one in \cite{si}.)
The following implies that measurable sets where the map is Lyapunov stable
are negligible with respect to any expansive measure (invariant or not).

\begin{prop}
\label{maluco}
If a measurable map of a separable metric space is
Lyapunov stable on a measurable set $A$, then $A$ has measure zero with respect to any
expansive measure.
\end{prop}

\noindent{{\bf Proof.}}
Fix a measurable map $f: X\to X$ of a separable metric space $X$,
an expansive measure $\mu$ and $\Delta>0$.
Since $\mu$ is regular there is a closed subset $C\subset A$
such that
$$
\mu(A\setminus C)\leq \Delta.
$$
Let us compute $\mu(C)$.

Fix an expansivity constant $\epsilon$ of $\mu$
(c.f. Definition \ref{sensitivity-constant}).
Since $f$ is Lyapunov stable on $A$ and $C\subset A$
for every $x\in C$ there is a neighborhood $U(x)$ such that
\begin{equation}
\label{fedorenko}
d(f^n(x),f^n(y))< \epsilon
\quad\quad\forall n\in \mathbb{N}, \forall y\in U(x)\cap C.
\end{equation}
On the other hand, $C$ is separable (since $X$ is) and so
Lindelof with the induced topology.
Consequently, the open covering $\{U(x)\cap C:x\in C\}$
of $C$ admits a countable subcovering $\{U(x_i)\cap C:i\in \mathbb{N}\}$.
Then,
\begin{equation}
\label{balurdo}
\mu(C)\leq\sum_{i\in \mathbb{N}}\mu\left(U(x_i)\cap C\right).
\end{equation}
Now fix $i\in \mathbb{N}$.
Applying (\ref{fedorenko}) to $x=x_i$ we obtain
$
U(x_i)\cap C\subset \Phi_\epsilon(x_i)
$
and then
$
\mu\left(U(x_i)\cap C\right)\leq \mu(\Phi_\epsilon(x))=0
$
since $\epsilon$ is an expansivity constant.
As $i$ is arbitrary we obtain $\mu(C)=0$ by (\ref{balurdo}).

To finish we observe that
$$
\mu(A)=\mu(A\setminus C)+\mu(C)=\mu(A\setminus C)\leq \Delta
$$
and so $\mu(A)=0$ since $\Delta$ is arbitrary.
This ends the proof.
\opensquare

\vspace{5pt}

From these propositions we obtain the following corollary.
Recall that
the {\em recurrent set} of $f: X\to X$ is
defined by $R(f)=\{x\in X:x\in \omega_f(x)\}$,
where
$
\omega_f(x)=\left\{y\in X:y=\lim_{k\to\infty}f^{n_k}(x)\mbox{ for some sequence }
n_k\to\infty\right\}.
$
\begin{clly}
\label{maps-zero-entropy-0}
A measurable map of a separable metric space which either
has countably many stable classes or
is Lyapunov stable on its recurrent set
has no expansive invariant measures.
\end{clly}

\noindent{{\bf Proof.}}
First consider the case when there are countably many stable classes.
Suppose by contradiction that there exists an expansive invariant measure.
Since the collection of stable classes is a partition of the space
it would follow from Proposition \ref{asymptotic-pair}
that the space has measure zero which is absurd.

Now consider the case when the map $f$ is Lyapunov stable on $R(f)$. Again
suppose by contradiction that there is an expansive invariant measure $\mu$.
Denote by $supp(\mu)$ the support of $\mu$.
Since $\mu$ is invariant we have $supp(\mu)\subset R(f)$
by Poincare recurrence. However, since $f$ is Lyapunov stable on $R(f)$ we obtain
$\mu(R(f))=0$ from Proposition \ref{maluco} so
$\mu(supp(\mu)))=\mu(R(f))=0$ which is absurd.
This proves the result.
\opensquare

\vspace{5pt}

\section{Applications}
\label{sec-entropy}

\noindent
We start this section by proving that positive entropy implies expansiveness
among ergodic invariant measures for continuous maps on compact metric spaces.
Afterward we include some short applications.

Recall that a Borel probability measure $\mu$ in $X$ is {\em ergodic} with respect to
a map if every invariant measurable set has measure zero or one.
For the notion of entropy $h_\mu(f)$ of an invariant measure $\mu$ see \cite{w}.

\begin{thm}
\label{thA}
Every ergodic invariant measure with positive entropy of a
continuous map on a compact metric space is expansive.
\end{thm}

\noindent{{\bf Proof.}}
Consider an ergodic invariant measure $\mu$ with positive entropy $h_\mu(f)>0$
of a continuous map $f$ on a compact metrix space $X$.
Fix $\delta>0$ and define
$$
X_\delta=\{x\in X:\mu(\Phi_\delta(x))=0\}.
$$
By Lemma \ref{forall} we are left to prove that there is 
$\delta>0$ such that $\mu(X_\delta)=1$.

Fix $x\in X$. It follows from the definition of $\Phi_\delta(x)$ that
$
\Phi_\delta(x)\subset f^{-1}(\Phi_\delta(f(x)))
$
so
$$
\mu(\Phi_\delta(x))\leq \mu(\Phi_\delta(f(x)))
$$
since $\mu$ is invariant.
Then, $\mu(\phi_\delta(x))=0$
whenever $x\in f^{-1}(X_\delta)$ yielding
$$
f^{-1}(X_\delta)\subset X_\delta.
$$
Denote by $A\Delta B$ the symmetric difference of the sets
$A,B$. Since $\mu(f^{-1}(X_\delta))=\mu(X_\delta)$ the above inequality implies that
$X_\delta$ is {\em essentially invariant}, i.e.,
$\mu(f^{-1}(X_\delta)\Delta X_\delta)=0$.
Since $\mu$ is ergodic we conclude that $\mu(X_\delta)\in\{0,1\}$ for all $\delta>0$.
Then, we are left to prove that there is $\delta>0$ such that $\mu(X_\delta)>0$.
To find it we proceed as follows.

For all $\delta>0$ we define the map $\phi_\delta: X\to I\!\! R\cup\{\infty\}$,
$$
\phi_\delta(x)=\liminf_{n\to\infty}-\frac{\log\mu(B[x,n,\delta])}{n}
$$
where $B[x,n,\delta]=\bigcap_{i=0}^{n-1}f^{-i}(B[f^i(x),\delta])$
and $B[\cdot,\cdot]$ stands for the closed ball operation.
Define $h=\frac{h_\mu(f)}{2}$
(thus $h>0$) and
$$
X^m=\left\{x\in X:\phi_{\frac{1}{m}}(x)>h\right\},
\quad\quad\forall m\in I\!\! N^+.
$$
Notice that $\phi_\delta(x)\geq \phi_{\delta'}(x)$ whenever
$0<\delta<\delta'$.
From this it follows that $X^m\subset X^{m'}$ for $m\leq m'$ and further
$$
\left\{
x\in X:\sup_{\delta>0}\phi_\delta(x)=h_\mu(f)
\right\}
\subset\bigcup_{m\in I\!\! N^+}X^m.
$$
Then,
$$
\mu
\left(
\left\{
x\in X:\sup_{\delta>0}\phi_\delta(x)=h_\mu(f)
\right\}
\right)\leq \lim_{m\to\infty}\mu(X^m).
$$
On the other hand, $\mu$ is nonatomic
since it is ergodic invariant with positive entropy.
So, the Brin-Katok Theorem \cite{bk} implies
$$
\mu
\left(
\left\{
x\in X:\sup_{\delta>0}\phi_\delta(x)=h_\mu(f)
\right\}
\right)
=1.
$$
Then,
$$
\lim_{m\to \infty}\mu(X^m)=1.
$$
Consequently, we can fix
$m\in I\!\! N^+$ such that
$$
\mu(X^m)>0.
$$
We shall prove that $\delta=\frac{1}{m}$ works.

Let us take $x\in X^m$.
It follows that $\mu(B[x,n,\delta])<e^{-hn}$ for all $n$ large.
Since $h>0$ we conclude that
$$
\lim_{n\to\infty}\mu(B[x,n,\delta])=0.
$$
But it follows from the definition of $\Phi_\delta(x)$
that
$$
\Phi_\delta(x)= \bigcap_{n=0}^\infty B[x,n,\delta].
$$
In addition, $B[x,n',\delta]\subset B[x,n,\delta]$ whenever $n\leq n'$ therefore
$$
\mu(\Phi_\delta(x))= \lim_{n\to\infty}\mu(B[x,n,\delta])=0.
$$
This proves $x\in X_\delta$.
As $x\in X^m$ is arbitrary we obtain $X^m\subset X_\delta$ whence
$$
0<\mu(X^m)\leq \mu(X_\delta)
$$
and the proof follows.
\opensquare

\vspace{5pt}

The converse of the above theorem is false, i.e., an
expansive measure may have zero entropy even in the ergodic invariant case.
A counterexample is as follows.

\begin{ex}
\label{circle-interval}
There are continuous maps in the circle exhibiting
ergodic invariant measures with zero entropy which, however, are expansive.
\end{ex}

\noindent{{\bf Proof.}}
Recall that a {\em Denjoy map} is a nontransitive homeomorphism with irrational rotation number of the circle $S^1$
(c.f. \cite{hk}).
Since all circle homeomorphisms have zero topological entropy it remains to prove
that every Denjoy map $h$ exhibits expansive measures.
As is well-known $h$ is uniquely ergodic and the support of its unique
invariant measure $\mu$ is a {\em minimal set}, i.e., a set which is minimal with respect to the property of being compact invariant.
We shall prove that this measure is expansive.
Denote by $E$ the support of $\mu$.
It is well known that $E$ is a Cantor set.
Let $\alpha$ be half of the length of the biggest interval $I$ in the complement
$S^1-E$ of $E$ and take $0<\delta<\alpha/2$.
Fix $x\in S^1$ and denote by $Int(\cdot)$ the interior operation.
We claim that $Int(\Phi_\delta(x))\cap E=\emptyset$.
Otherwise, there is some $z\in Int(\Phi_\delta(x))\cap E$.
Pick $w\in \partial I$ (thus $w\in E$).
Since $E$ is minimal there is a sequence $n_k\to\infty$ such that
$h^{-n_k}(w)\to z$. Since $\mu$ is a finite measure, the interval sequence
$\{h^{-n}(I):n\in I\!\! N\}$ is disjoint, we have that
the length of the intervals $h^{-n_k}(I)\to 0$ as $k\to\infty$.
It turns out that there is some integer $k$ such that
$h^{-n_k}(I)\subset \Phi_\delta(x)$.
From this and the fact that $h(\Phi_\delta(x))\subset \Phi_\delta(h(x))$
one sees that $I\subset B[h^{n_k}(x),\delta]$ which is clearly absurd because
the length of $I$ is greather than $\alpha>2\delta$.
This contradiction proves the claim.
Since $\Phi_\delta(x)$ is either a closed interval or $\{x\}$ the
claim implies that $\Phi_\delta(x)\cap E=\Phi_\delta(x)\cap E$
consists of at most two points.
Since $\mu$ is clearly nonatomic we conclude that $\mu(\Phi_\delta(x))=0$.
Since $x\in S^1$ is arbitrary we are done.
\opensquare

\vspace{5pt}

\begin{rk}
Altogether the above proof and \cite{m} characterize
the Denjoy maps as those circle homeomorphisms exhibiting expansive measures.
\end{rk}

A first application of Theorem \ref{thA} is as follows.

\begin{thm}
\label{stable-set-0}
The stable classes of a continuous map of a compact metric space
have measure zero with respect to any ergodic invariant measure with positive entropy.
\end{thm}

\noindent{{\bf Proof.}}
In fact, since these measures are expansive by Theorem \ref{thA}
we obtain the result from Proposition \ref{asymptotic-pair}.
\opensquare

\vspace{5pt}

We can also use Theorem \ref{thA} to compute the topological entropy of certain
continuous maps on compact metric spaces (for the related concepts see \cite{akm} or \cite{w}).
As a motivation let us mention the known facts that both
{\em transitive} continuous maps with countably many stable classes on compact metric spaces
and continuous maps of the {\em interval} or the {\em circle} which are Lyapunov stable on their
recurrent sets
have zero topological entropy
(see Corollary 2.3 p. 263 in \cite{hy},
\cite{fs}, Theorem B in \cite{si} and \cite{z}).
Indeed we improve these result in the following way.

\begin{thm}
\label{maps-zero-entropy}
A continuous map of a compact metric space which either
has countably many stable classes or
is Lyapunov stable on its recurrent set
has zero topological entropy.
\end{thm}

\noindent{{\bf Proof.}}
If the topological entropy were not zero
the variational principle \cite{w} would imply the existence of
ergodic invariant measures with positive entropy.
But by Theorem \ref{thA} these measures are expansive against Corollary \ref{maps-zero-entropy-0}.
\opensquare

\vspace{5pt}

\begin{ex}
An example satisfying the first part of Theorem \ref{maps-zero-entropy} is
the classical pole North-South diffeomorphism on spheres.
In fact, the only stable sets of this diffeomorphism
are the stable sets of the
poles.
The Morse-Smale diffeomorphisms \cite{hk}
are basic examples where these hypotheses are fulfilled.
\end{ex}

Now we use expansive measures to study the chaoticity
in the sense of Li and Yorke \cite{liy}.
Recall that if $\delta\geq 0$ a {\em $\delta$-scrambled set} of $f:X\to X$
is a subset $S\subset X$ satisfying
\begin{equation}
 \label{delta-scrambled}
\liminf_{n\to\infty}d(f^n(x),f^n(y))=0\quad
\mbox{ and }
\quad
\limsup_{n\to\infty}d(f^n(x),f^n(y))>\delta
\end{equation}
for all different points $x,y\in S$.
Recalling that a {\em Polish space} is a complete separable metric space we obtain the following result.

\begin{thm}
 \label{cadre->li}
A continuous map of a Polish space carrying an uncountable $\delta$-scrambled set for some $\delta>0$
also carries expansive measures.
\end{thm}

\noindent{{\bf Proof.}}
Let $X$ a Polish space and $f: X\to X$ be a continuous map carrying
an uncountable $\delta$-scrambled set for some $\delta>0$.
Then, by Theorem 16 in \cite{bhs}, there is a {\em closed} uncountable $\delta$-scrambled set $S$.
As $S$ is closed and $X$ is Polish we have that $S$ is also
a Polish space with respect to the induced metric.
As $S$ is uncountable we have from \cite{prv} that there is a nonatomic Borel probability measure
$\nu$ in $S$.
Let $\mu$ be the Borel probability induced by $\nu$ in $X$, i.e.,
$\mu(A)=\nu(A\cap S)$ for all Borelian $A\subset X$.
We shall prove that this measure is expansive.
If $x\in S$ and $y\in \Phi_{\frac{\delta}{2}}(x)\cap S$
we have that $x,y\in S$ and
$d(f^n(x),f^n(y))\leq \frac{\delta}{2}$ for all $n\in \mathbb{N}$
therefore $x=y$ by the second inequality in (\ref{delta-scrambled}).
We conclude that $\Phi_{\frac{\delta}{2}}(x)\cap S=\{x\}$ for all $x\in S$.
As $\nu$ is nonatomic we obtain $\mu(\Phi_{\frac{\delta}{2}}(x))=\nu(\Phi_{\frac{\delta}{2}}(x)\cap S)=\nu(\{x\})=0$
for all $x\in S$.
On  other hand, it is clear that every open set which does not intersect $S$ has $\mu$-measure $0$
so $\mu$ is supported in the closure of $S$. As $S$ is closed we obtain that
$\mu$ is supported on $S$.
We conclude that
$\mu(\Phi_{\frac{\delta}{2}}(x))=0$ for $\mu$-a.e. $x\in X$,
so, $\mu$ is expansive by Lemma \ref{forall}.
\opensquare

\vspace{5pt}

Now recall that a continuous map is {\em Li-Yorke chaotic} if
it has an uncountable $0$-scrambled set.

Until the end of this section $M$ will denote either the interval $I=[0,1]$ or the unit circle $S^1$.

\begin{clly}
Every Li-Yorke chaotic map in $M$ carries expansive measures.
\end{clly}

\noindent{{\bf Proof.}}
Theorem in p. 260 of \cite{d} together with
theorems A and B in \cite{ku} imply that every Li-Yorke chaotic map in $I$
or $S^1$ has an uncountable $\delta$-scrambled set for some $\delta>0$.
Then, we obtain the result from Theorem \ref{cadre->li}.
\opensquare

\vspace{5pt}

It follows from Example \ref{circle-interval} that there are continuous maps with zero topological entropy
in the circle exhibiting expansive {\em invariant} measures.
This leads to the question whether the same result is true on compact intervals.
The following consequence of the above corollary gives a partial positive answer for this question.

\begin{ex}
There are continuous maps with zero topological entropy in the interval
carrying expansive measures.
\end{ex}

Indeed, by \cite{j} there is a continuous map of the interval, with zero topological
entropy, exhibiting a $\delta$-scrambled set of positive Lebesgue measures for some $\delta>0$.
Since sets with positive Lebesgue measure are uncountable we obtain an expansive measure
from Theorem \ref{cadre->li}.

Another interesting example is the one below.

\begin{ex}
 \label{tent}
The Lebesgue measure is an ergodic invariant measure with positive entropy of
the tent map $f(x)=1-|2x-1|$ in $I$. Therefore, this measure is expansive by Theorem \ref{thA}.
\end{ex}

It follows from this example that there are continuous maps in $I$
carrying expansive measures $\mu$ {\em with full support} (i.e. $supp(\mu)=I$).
These maps
also exist in $S^1$
(e.g. an expanding map).
Now, we prove that Li-Yorke and positive topological entropy are
equivalent properties among these maps in $I$.
But previously we need a result based
on the following well-known definition.

A {\em wandering interval}
of a map $f: M\to M$ is an interval $J\subset M$
such that $f^n(J)\cap f^m(J)=\emptyset$ for all different integers $n,m\in \mathbb{N}$
and no point in $J$ belongs to the stable set of some periodic point.

\begin{lemma}
\label{non-I}
If $f: M\to M$ is continuous,
then every wandering interval has measure zero with respect to every expansive measure.
\end{lemma}

\noindent{{\bf Proof.}}
Let $J$ a wandering interval and $\mu$ be an expansive measure with expansivity constant $\epsilon$
(c.f. Definition \ref{sensitivity-constant}).
To prove $\mu(J)=0$ it suffices to prove
$Int(J)\cap supp(\mu)=0$ since $\mu$ is nonatomic.
As $J$ is a wandering interval one has
$\lim_{n\to\infty}|f^n(J)|=0$,
where $|\cdot |$ denotes the length operation.
From this there is a positive integer $n_0$ satisfying
\begin{equation}
\label{mao}
|f^n(J)|<\epsilon,
\quad\quad\forall n\geq n_0.
\end{equation}
Now, take $x\in Int(J)$.
Since $f$ is clearly uniformly continuous and $n_0$ is fixed we can select
$\delta>0$ such that $B[x,\delta]\subset Int(J)$
and $|f^n(B[x,\delta])|<\epsilon$ for
$0\leq n\leq n_0$.
This together with
(\ref{mao}) implies
$|f^n(x)-f^n(y)|<\epsilon$ for all $n\in\mathbb{N}$ therefore
$B[x,\delta]\subset \Phi_\epsilon(x)$
so $\mu(B[x,\delta])=0$ since $\epsilon$ is an expansivity constant.
Thus $x\not\in supp(\mu)$ and we are done.
\opensquare

\vspace{5pt}

From this we obtain the following corollary.

\begin{clly}
A continuous map with expansive measures of the circle or the interval has no wandering intervals.
Consequently, a continuous map of the interval carrying expansive measures
with full support is Li-Yorke chaotic if and only if
it has positive topological entropy.
\end{clly}

\noindent{{\bf Proof.}}
The first part is a direct consequence Lemma \ref{mao}
while, the second, follows from the first
since a continuous interval map without wandering intervals is Li-Yorke chaotic
if and only if it has positive topological entropy \cite{smi}.
\opensquare


\vspace{5pt}













\section*{References}
\begin{thebibliography}{00}


\bibitem{akm}
Adler, R., L., Konheim, A., G., McAndrew, M., H.,
{\em Topological entropy}, Trans. Amer. Math. Soc. 114 (1965), 309--319.





\bibitem{bhr}
Blanchard, F., Host, B., Ruette, S.,
{\em Asymptotic pairs in positive-entropy systems},
Ergodic Theory Dynam. Systems 22 (2002), no. 3, 671--686.









\bibitem{bhs}
Blanchard, F., Huang, W., Snoha, L.,
{\em Topological size of scrambled sets}.
Colloq. Math. 110 (2008), no. 2, 293--361.







\bibitem{b}
Bowen, R.,
{\em Entropy-expansive maps},
Trans. Amer. Math. Soc.  164  (1972), 323--331.




\bibitem{bk}
Brin, M., Katok, A.,
{\em On local entropy},
Geometric dynamics (Rio de Janeiro, 1981), 30--38, Lecture Notes in Math., 1007, Springer, Berlin, 1983.









\bibitem{cj}
Cadre, B., Jacob, P.,
{\em On pairwise sensitivity},
J. Math. Anal. Appl. 309 (2005), no. 1, 375--382.




\bibitem{d}
Du, B-S.
{\em Every chaotic interval map has a scrambled set in the recurrent set},
Bull. Austral. Math. Soc. 39 (1989), no. 2, 259--264.








\bibitem{fs}
Fedorenko, V., V., Smital, J.,
{\em Maps of the interval Ljapunov stable on the set of nonwandering points},
Acta Math. Univ. Comenian. (N.S.) 60 (1991), no. 1, 11--14.





\bibitem{hk}
Hasselblatt, B., Katok, A.,
{\em Introduction to the modern theory of dynamical systems. With a supplementary
chapter by Katok and Leonardo Mendoza},
Encyclopedia of Mathematics and its Applications, 54. Cambridge University Press, Cambridge, 1995.




\bibitem{h}
Huang, W.,
{\em Stable sets and $\epsilon$-stable sets in positive-entropy systems},
Comm. Math. Phys. 279 (2008), no. 2, 535--557.



\bibitem{hy}
Huang, W., Ye, X.,
{\em Devaney's chaos or 2-scattering implies Li-Yorke's chaos},
Topology Appl. 117 (2002), no. 3, 259--272.






\bibitem{j}
Jimenez Lopez, V.,
{\em Large chaos in smooth functions of zero topological entropy},
Bull. Austral. Math. Soc. 46 (1992), no. 2, 271--285.












\bibitem{ku}
Kuchta, M.,
{\em Characterization of chaos for continuous maps of the circle},
Comment. Math. Univ. Carolin. 31 (1990), no. 2, 383--390.





\bibitem{liy}
Li, T.,Y., Yorke, J., A.,
{\em Period three implies chaos},
Amer. Math. Monthly 82 (1975), no. 10, 985--992.




\bibitem{m1}
Morales, C., A.,
{\em A generalization of expansivity},
Discrete Contin. Dyn. Syst. (to appear).


\bibitem{m}
Morales, C., A.,
{\em Measure-expansive systems},
Preprint IMPA (2011), D083 / 2011.













\bibitem{prv}
Parthasarathy, K., R., Ranga Rao, R., Varadhan, S., R., S.,
{\em On the category of indecomposable distributions on topological groups},
Trans. Amer. Math. Soc. 102 (1962), 200--217.





\bibitem{si}
Sindelarova, P.,
{\em A counterexample to a statement concerning Lyapunov stability}
Acta Math. Univ. Comenianae 70 (2001), 265--268.




\bibitem{smi}
Smital, J.,
{\em Chaotic functions with zero topological entropy},
Trans. Amer. Math. Soc. 297 (1986), no. 1, 269--282.


\bibitem{u}
Utz, W., R.,
{\em Unstable homeomorphisms},
Proc. Amer. Math. Soc.  1 (1950), 769--774.




\bibitem{w}
Walters, P.,
{\em An introduction to ergodic theory},
Graduate Texts in Mathematics, 79. Springer-Verlag, New York-Berlin, 1982.






\bibitem{z}
Zhou, Z., L.,
{\em Some equivalent conditions for self-mappings of a circle},
(Chinese) Chinese Ann. Math. Ser. A 12 (1991), suppl., 22--27.




\end{thebibliography}
\end{document}